\RequirePackage{fix-cm}
\documentclass[smallextended]{svjour3}       
\smartqed  
\usepackage{graphicx} 
\usepackage{amsmath,amssymb,amsfonts,amstext}
\usepackage{latexsym,graphics,epsf,epsfig,psfrag}
\usepackage{color}
\usepackage{thmtools}
\usepackage[linesnumbered, ruled]{algorithm2e}
\SetKwRepeat{Do}{do}{while}%
\usepackage[noend]{algpseudocode}
\usepackage{dsfont}
\usepackage[colorlinks,linkcolor=black,citecolor=black,urlcolor=black]{hyperref}
\usepackage{enumitem}
\usepackage{balance}

\newtheorem{thm}{Theorem}
\newtheorem{coro}{Corollary}
\newtheorem{assum}{Assumption}

\usepackage{xspace}

\newcommand{\ie}{{i.e.}\xspace}

\newcommand{\ab}{\mathbf{a}}

\newcommand{\xb}{\mathbf{x}}
\newcommand{\yb}{\mathbf{y}}

\newcommand{\ub}{\mathbf{u}}

\newcommand{\RR}{\mathds{R}}

\newcommand{\dom}{\mathop{\mathrm{dom}}}

\newcommand{\argmin}{\mathop{\mathrm{argmin}}}
\newcommand{\inner}[2]{\langle #1, #2 \rangle}

\begin{document}

\title{A Simple Convergence Analysis of Bregman Proximal Gradient Algorithm \thanks{The work of Y. Zhou and Y. Liang was supported by the National Science Foundation under Grant CCF-1704169 and ECCS-1609916. The work of L. Shen was supported in part  by the National Science Foundation under grant DMS-1115523 and DMS-1522332, and by the National Research Council.}
}


\author{Yi Zhou  \and
        Yingbin Liang \and
        Lixin Shen 
}


\institute{Yi Zhou \at
              Ohio State University, Columbus, OH, USA. \\
              Department of Electrical and Computer Engineering. \\
              \email{zhou.1172@osu.edu}          
           \and
             Yingbin Liang \at
             Ohio State University, Columbus, OH, USA. \\
             Department of Electrical and Computer Engineering.\\
             \email{liang.889@osu.edu}     
           \and 
           Lixin Shen \at
           Syracuse University, Syracuse, NY, USA. \\
           Department of Mathematics.  \\
           \email{lshen03@syr.edu}   
}

\date{Received: date / Accepted: date}

\maketitle

\begin{abstract}
In this paper, we provide a simple convergence analysis of proximal gradient algorithm with Bregman distance, which provides a tighter bound than existing result. In particular, for the problem of minimizing a class of convex objective functions, we show that proximal gradient algorithm with Bregman distance can be viewed as proximal point algorithm that incorporates another Bregman distance. Consequently, the convergence result of the proximal gradient algorithm with Bregman distance follows directly from that of the proximal point algorithm with Bregman distance, and this leads to a simpler convergence analysis with a tighter convergence bound than existing ones. We further propose and analyze the backtracking line search variant of the proximal gradient algorithm with Bregman distance. Simulation results show that the line search method significantly improves the convergence performance of the algorithm.

\end{abstract}

\section{Introduction}\label{intro}
Proximal algorithms have been extensively studied in optimization theory, as they are efficient solvers to problems that involve non-smoothness and have a fast convergence rate. The proximal algorithms have been widely applied to solve practical problems including image processing, e.g., \cite{FISTA,Micchelli-Shen-Xu:IP-11}, distributed statistical learning, e.g., \cite{distributed-admm}, and low rank matrix minimization, e.g., \cite{low-rank}.

Consider the following  optimization problem:
\begin{flalign}\label{eq:problem}
\textbf{(P1)}\quad \min_{\xb \in C}~r(\xb),
\end{flalign}
where $r:\RR^n \to (-\infty, +\infty]$ is a proper, lower-semicontinuous convex function, and $C$ is a closed convex set in $\RR^n$. The well known proximal point algorithm (PPA) for solving \textbf{(P1)} was introduced initially by Martinet in \cite{Martinet-PPA}. The algorithm generates a sequence $\{\xb_k\}$ via the following iterative step:
\begin{flalign}\label{eq:PPA}
\text{(PPA-$\mathcal{E}$)}~~ \xb_{k+1}=\underset{\xb \in C}{\argmin}\left\{r(\xb) + \frac{1}{2\lambda_{k}}\|\xb-\xb_{k}\|_2^{2}\right\},
\end{flalign}
where $\lambda_k>0$ corresponds to the step size at $k$-th iteration. We refer to the algorithm as PPA-$\mathcal{E}$ for the choice of the Euclidean distance (\ie, the $\|\cdot\|_2^2$ term).
This algorithm can be interpreted as applying the gradient descent method on the Moreau envelope of $r$, \ie, a smoothed version of the objective function \cite{Nesterov-smooth,smooth-Marc}.
It was shown in \cite{RT-PPA,Eckstein-DR} that the sequence $\{\xb_k\}$ generated by PPA-$\mathcal{E}$ converges to a solution of \textbf{(P1)} with a proper choice of the step size sequence $\{\lambda_k\}$, and the rate of convergence was characterized in \cite{Guler-PPA}.

A natural generalization of PPA-$\mathcal{E}$ is to replace the Euclidean distance with a more general distance-like term. In existing literature, various choices of distance have been proposed, e.g., \cite{Censor-GPPA,Entro-Marc,Marc-proximal-like,Tebolle_mirror}. Among them, a popular choice is the Bregman distance \cite{GPPA-Chen-Marc,BPPA-Eckstein,Tebolle_mirror}.
As a consequence, we obtain the following iterative step of proximal point algorithm with Bregman distance:
\begin{flalign}\label{eq:GPPA}
\text{(PPA-$\mathcal{B}$)}\quad \xb_{k+1}=\underset{\xb \in C}{\argmin}\left\{r(\xb) + \frac{1}{\lambda_{k}}D_h(\xb,\xb_{k})\right\}.
\end{flalign}
Here, $D_h(\xb,\xb_k)$ corresponds to the Bregman distance between the points $\xb$ and $\xb_k$, and is based on a continuously differentiable strictly convex function $h$. We refer to the algorithm as PPA-$\mathcal{B}$ for the choice of the Bregman distance, which is formally defined in Definition \ref{def:Bregman distance class}. The convergence rate of PPA-$\mathcal{B}$ has been characterized in \cite{GPPA-Chen-Marc,Marc-proximal-like}, and we refer to \cite{GPPA-Marc} for a comprehensive discussion on the PPA with different choices of distance metrics.


A generalized optimization problem of \textbf{(P1)} is the following composite objective minimization problem:
\begin{flalign}\label{eq:problem f+g}
\textbf{(P2)}\quad \min_{\xb \in C} ~\{F(\xb):= f(\xb)+g(\xb)\},
\end{flalign}
where $f$ is usually a differentiable and convex loss function that corresponds to the data fitting part, and $g$ is a possibly non-smooth regularizer that promotes structures such as sparsity, low-rankness, etc, to the solution of the problem. This composite objective minimization problem generalizes many applications in machine learning, image processing, detection, etc.
As an extension of the PPA, splitting algorithms are proposed for solving the composite objective minimization problem in \textbf{(P2)} \cite{Eckstein-DR,splitting-algo}.
In particular, the proximal gradient algorithm (PGA) with Euclidean distance has been developed in \cite{goldstein-PGA} to solve \textbf{(P2)} efficiently, and the iterative step is given by
\begin{flalign}\label{eq:PGA}
\text{(PGA-$\mathcal{E}$)}: \quad \xb_{k+1}=\underset{\xb \in C}{\argmin}\left\{g(\xb) + \langle \xb, \nabla f(\xb_{k}) \rangle +\frac{1}{2\gamma_k}\|\xb - \xb_{k}\|_{2}^{2}\right\},
\end{flalign}
where we refer to the algorithm as PGA-$\mathcal{E}$ for the choice of Euclidean distance.
It has been shown that the sequence of function value residual generated by PGA-$\mathcal{E}$ has a convergence rate of $\mathcal{O}(1/k)$\footnote{Here, $f(n)=\mathcal{O}(g(n))$ denotes that $|f(n)|\le \xi |g(n)|$ for all $n>N$, where $\xi$ is a constant and $N$ is a positive integer.} \cite{FISTA}, and the rate can be further improved to be $\mathcal{O}(1/k^{2})$ via Nesterov's acceleration technique. Inspired by the way of generalizing PPA-$\mathcal{E}$ to PPA-$\mathcal{B}$, PGA-$\mathcal{E}$ can also be generalized by replacing the Euclidean distance with the Bregman distance, and correspondingly, the iterative step is given by
\begin{flalign}\label{eq:BPGA}
\text{(PGA-$\mathcal{B}$)}: \quad \xb_{k+1}=\underset{\xb \in C}{\argmin}\left\{g(\xb) + \langle \xb, \nabla f(\xb_{k}) \rangle + \frac{1}{\gamma_k}D_{h}(\xb,\xb_{k})\right\},
\end{flalign}
where we refer to the algorithm as PGA-$\mathcal{B}$ for the choice of Bregman distance.
Under Lipschitz continuity of $\nabla f$, it has been shown in \cite{Tseng-BPGA} that the convergence rate of PGA-$\mathcal{B}$ is $\mathcal{O}(1/k)$.

It is clear that PGA-$\mathcal{E}$ and PGA-$\mathcal{B}$ are respectively generalizations of PPA-$\mathcal{E}$ and PPA-$\mathcal{B}$, because they coincide when the function $f$ in \textbf{(P2)} is a constant function. Thus, in existing literature, the analysis of PGA is developed by their own as in \cite{FISTA,Tseng-BPGA} without resorting to existing analysis of PPA. More recently, a concurrent work \cite{bolte_BPGA} to this paper interprets PGA-$\mathcal{B}$ as the composition of mirror descent method and PPA-$\mathcal{B}$.
In contrast to this viewpoint, this paper shows that PGA-$\mathcal{B}$ can, in fact, be viewed as PPA-$\mathcal{B}$ under a proper choice of the Bregman distance. Consequently, the analysis of PGA-$\mathcal{B}$ can be mapped to that of PPA-$\mathcal{B}$, resulting in a much simpler analysis. We note that the initial version of this paper \cite{Yi2015} was posted on arXiv in March, 2015, which already independently developed the aforementioned main result.

We summarize our main contributions as follows.
In this paper, we show that PGA-$\mathcal{B}$ can be viewed as PPA-$\mathcal{B}$ with a special choice of Bregman distance, and thus the convergence results of PGA-$\mathcal{B}$ inherit existing convergence results of PPA-$\mathcal{B}$. Following this viewpoint, we obtain a tighter bound of the convergence rate of the function value residual, and our result avoids involving the symmetry coefficient in \cite{bolte_BPGA}. Lastly, we propose a line search variant of PGA-$\mathcal{B}$ and characterize its convergence rate.


The rest of the paper is organized as follows. In Section \S\ref{sec: distance}, we recall the definition of Bregman distance, unify PGA-$\mathcal{B}$ as a special case of PPA-$\mathcal{B}$ and discuss its convergence results. In Section \S\ref{sec: line search}, we propose a line search variant of PGA-$\mathcal{B}$ and characterize its convergence rate. In Section \S \ref{sec: exp}, we compare the convergence behavior between PGA-$\mathcal{B}$ and its line search variant via numerical experiments.
Finally in Section \S \ref{sec:conclude}, we conclude our paper with a few remarks on our results.

\section{Unifying PGA-$\mathcal{B}$ as PPA-$\mathcal{B}$}\label{sec: distance}
\subsection{Preliminaries on Bregman Distance}
We first recall the definition of the Bregman distance \cite{Bregman1967}, see also \cite{Pierro1986,GPPA-Chen-Marc}. Throughout, the interior of a set $C\subset \RR^n$ is denoted as $\mathrm{int}C$.
\begin{definition}[Bregman Distance $\mathcal{B}$] \label{def:Bregman distance class}
	Let $h: \RR^n \to (-\infty, +\infty]$ be a function with $\dom h = C$, and satisfies:
	\begin{enumerate}[topsep=0pt,noitemsep]
		\item[(a)] $h$ is continuously differentiable on $\mathrm{int}C$;
		\item[(b)] $h$ is strictly convex on $C$.
	\end{enumerate}
	Then, the Bregman distance $D_{h}: C \times \mathrm{int}C \to \RR_+$ associated with function $h$ is defined as, for all $\xb \in C$ and $\yb \in \mathrm{int}C$,
	\begin{flalign}\label{eq:Bregman distance}
	D_{h}(\xb,\yb) = h(\xb) - h(\yb) - \langle \xb-\yb, \nabla h(\yb)\rangle.
	\end{flalign}
\end{definition}
We denote $\mathcal{B}$ as the class of all Bregman distances.
Clearly, the Bregman distance $D_h$ is defined as the residual of the first order Taylor expansion of function $h$. In general, the Bregman distance is asymmetric with respect to the two arguments. On the other hand, the convexity of function $h$ implies the non-negativity of the Bregman distance, making it behaves like a metric. Moreover, the following properties are direct consequences of \eqref{eq:Bregman distance}: For any $\ub\in C, \xb, \yb \in \mathrm{int}C$ and any $D_h, D_{h'} \in \mathcal{B}$,
\begin{flalign}
&D_{h}(\ub,\xb) + D_{h}(\xb,\yb) - D_{h}(\ub,\yb) = \langle \nabla h(\yb) - \nabla h(\xb), \ub - \xb \rangle. \label{eq:three point}\\
&D_{h}(\ub,\xb) \pm D_{h'}(\ub,\xb) = D_{h\pm h'}(\ub,\xb). \label{eq:linearity}
\end{flalign}
The property in \eqref{eq:three point} establishes a relationship among the Bregman distances of three points, and the property in \eqref{eq:linearity} shows the linearity of the Bregman distance with respect to the function $h$. In summary, Bregman distances are similar to metrics (but they can be asymmetric), and the following are several popular examples of Bregman distance.


\begin{example}(Euclidean Distance)\label{ex: Eucli}
	For $h: \RR^n \to \RR$ with $h(\xb) = \frac{1}{2}\|\xb\|_2^2$, $D_h(\xb, \yb) = \frac{1}{2}\|\xb-\yb\|_2^2$.
\end{example}
\begin{example}(KL Relative Entropy)\label{ex: KL}
	For $h: \RR^n_+ \to \RR$ with $h(\xb) = \sum_{j=1}^{n} x_j\log x_j - x_j$ (with the convention $0\log0 = 0$), $D_{h}(\xb, \yb) = \sum_{j=1}^{n} x_j \log \frac{x_j}{y_j} - x_j + y_j$.
\end{example}
\begin{example}(Burg's Entropy)\label{ex: burg}
	For $h: \RR^n_{++} \to \RR$ with $h(\xb) = -\sum_{i=1}^{n} \log x_i$, $D_h(\xb, \yb) = \sum_{j=1}^{n} \frac{x_j}{y_j} - \log \frac{x_j}{y_j} - 1 $.
\end{example}
Clearly, the Euclidean distance in Example \ref{ex: Eucli} is a special case of Bregman distance, and hence the proximal algorithms under Bregman distance naturally generalize the corresponding ones under Euclidean distance. The Bregman distances in Example \ref{ex: KL} and \ref{ex: burg} have a non-Euclidean structure. In particular, the Kullback-Liebler (KL) relative entropy is useful when the set $C$ is the simplex \cite{Tebolle_mirror}, and the Burg's Entropy is suitable for optimizing Poisson log-likelihood functions \cite{bolte_BPGA}.

\subsection{Connecting PGA-$\mathcal{B}$ to PPA-$\mathcal{B}$}\label{sec:BPGA}
Consider applying PGA-$\mathcal{B}$ to solve \textbf{(P2)}. The following standard assumptions are adopted regarding the functions $f,g,h$.
\begin{assum}\label{assum: compo}
	Regarding $f,g,h:\RR^n \to (-\infty, +\infty]$:
	\begin{enumerate}
		\item Functions $f,g$ are proper, lower semicontinuous and convex functions, $f$ is differentiable on $\mathrm{int}C$; $\dom f \supset C, \dom g \cap \mathrm{int}C \ne \emptyset$;
		\item $F^* := \inf_{\xb\in C} F(\xb) > - \infty$, and the solution set $\mathcal{X}^* := \{\xb: ~F(\xb) = F^*\}$ is non-empty;
		\item Function $h$ satisfies the properties in Definition \ref{def:Bregman distance class}.
	\end{enumerate}
\end{assum}

To simplify the analysis, we also assume that the iteration step of PGA-$\mathcal{B}$ is well defined, and refer to \cite{GPPA-Marc,bolte_BPGA,Censor-GPPA} for a detailed discussion. The following theorem establishes the main result that connects PGA-$\mathcal{B}$ with PPA-$\mathcal{B}$.
\begin{thm}\label{thm:BFBS}
	Assume the iteration steps of PGA-$\mathcal{B}$ are well defined. Then, the iteration step of PGA-$\mathcal{B}$ in \eqref{eq:BPGA} is equivalent to the following PPA-$\mathcal{B}$ step:
	\begin{flalign}\label{eq:thm-BFBS}
	\xb_{k+1}=\underset{\xb \in C}{\argmin}\left\{F(\xb)+D_{\ell_k}(\xb,\xb_{k})\right\},
	\end{flalign}
	where the function $\ell_k = \tfrac{1}{\gamma_k} h - f $.
\end{thm}
\begin{proof}
	By linearity in \eqref{eq:linearity} and the definition of Bregman distance in \eqref{eq:Bregman distance}, we obtain
	\begin{flalign}
	F(\xb)+D_{\ell_k}(\xb,\xb_{k}) &=f(\xb)+g(\xb)+\tfrac{1}{\gamma_k}D_{h}(\xb,\xb_{k})-D_{f}(\xb,\xb_{k}) \nonumber\\
	&=g(\xb) \!+\! \langle \xb, \nabla f(\xb_{k}) \rangle \!+\!\tfrac{1}{\gamma_k}D_{h}(\xb,\xb_{k})\!+\!f(\xb_k)\!-\!\langle \xb_k, \nabla f(\xb_{k}) \rangle.\nonumber
	\end{flalign}
	Thus, by ignoring the last two constant terms, the minimization problem of PGA-$\mathcal{B}$ is equivalent to \eqref{eq:thm-BFBS}, which is a PPA-$\mathcal{B}$ step with Bregman distance $D_{\ell_k}$. This completes the proof.
\end{proof}
Thus, PGA-$\mathcal{B}$ can be  mapped  exactly into the form of PPA-$\mathcal{B}$, and the form in \eqref{eq:thm-BFBS} provides a new insight of PGA---It is PPA with a special Bregman distance $D_{\ell_k}$. In particular, the $-f$ part of the function $\ell_k$ linearizes the objective function $f$, \ie, $f(\xb) + D_{-f}(\xb,\xb_k) = f(\xb_k) + \inner{\xb - \xb_k}{\nabla f(\xb_k)}$, which is a linear function. The linearizion simplifies the subproblem at each iteration, and leads to an update rule with closed form for a simple regularizer $g$.

We can further understand the class of objective functions that can be solved by PGA-$\mathcal{B}$ with a theoretical guarantee by leveraging the above equivalence viewpoint. In particular, to make the equivalent PPA-$\mathcal{B}$ step in \eqref{eq:thm-BFBS} be proper, the function $\ell_k$ in the Bregman distance should be convex and independent of the iteration $k$. Then, we are motivated to make the following assumption on the composite objective function.
\begin{assum}\label{assum: lipschitz}
	For the function $f$ in problem \textbf{(P2)}, there exists $\bar{\gamma} > 0$ and function $h$ in Definition \ref{def:Bregman distance class} such that for $\gamma_k = \gamma, k = 1, 2, \ldots$ with $0< \gamma < \bar{\gamma}$, the function $\ell_k := \ell = \frac{1}{\gamma}h - f$ is convex on $C$.
\end{assum}
Here, we consider the case $\gamma_k \equiv \gamma$, which corresponds to the choice of constant step size of PGA-$\mathcal{B}$. We further provide a backtracking line search rule for choosing the stepsize in Section \ref{sec: line search}.
Assumption \ref{assum: lipschitz} has also been considered in \cite{bolte_BPGA} to generalize the assumptions that
$\nabla f$ is Lipschitz continuous respect to certain norm, with respect to which function $h$ is strongly convex. In comparison, Assumption \ref{assum: lipschitz} does not require function $h$ and $\nabla f$ to satisfy these structures under certain norm. This generalization is useful, as some practical problems have objective functions that do not have norm structures. This point is further illustrated by Example \ref{ex: poisson}, which is presented after the convergence results. 

Our view point of PGA-$\mathcal{B}$ is very different from that developed in a concurrent independent work \cite{bolte_BPGA}. There, they view PGA-$\mathcal{B}$ as a mirror descent step composed with a PPA-$\mathcal{B}$ step, and develop a generalized descent lemma based on Assumption \ref{assum: lipschitz} to analyze the algorithm.  Our approach, however, is straightforward --- we simply map the PGA-$\mathcal{B}$ step exactly to a PPA-$\mathcal{B}$ step under Assumption \ref{assum: lipschitz}. This provides a unified view of PGA-$\mathcal{B}$ as a special case of PPA-$\mathcal{B}$, and consequently, the convergence analysis of PGA-$\mathcal{B}$ naturally follows from those of PPA-$\mathcal{B}$. In particular, Lemma 3.3 of \cite{GPPA-Chen-Marc} proposed the following properties of PPA-$\mathcal{B}$.
\begin{lemma}\cite[Lemma 3.3]{GPPA-Chen-Marc}\label{lemma: ppa-b}
	Consider the problem \textbf{(P1)} with optimal solution set $\mathcal{X}^*$. Let $\{\lambda_k\}$ be a sequence of positive numbers and denote $\sigma_k:=\sum_{l=1}^{k}\lambda_l$. Then the sequence $\{\xb_k\}$ generated by PPA-$\mathcal{B}$ given in \eqref{eq:GPPA} satisfy  
	\begin{align}
	&r(\xb_{k+1}) - r(\xb_{k}) \le - D_{h}(\xb_k, \xb_{k+1}),\label{eq: func_dec-ppa}\\
	& D_{h}(\xb^*, \xb_{k+1}) \le D_{h}(\xb^*, \xb_{k}), \quad \forall \xb^* \in \mathcal{X}^*, \label{eq: dist_dec-ppa}\\
	&r(\xb_{k}) - r(\ub) \le \frac{D_{h}(\ub,\xb_{0})}{k}, \quad\forall \ub \in C,\label{eq: global-ppa}
	\end{align}
\end{lemma}

By the connection between PPA-$\mathcal{B}$ and PGA-$\mathcal{B}$ that established in Theorem \ref{thm:BFBS}, we now identify $r = F, \lambda_k \equiv 1 (\sigma_k = k), h = \ell$ in Lemma \ref{lemma: ppa-b}, and directly obtain the following results on the iterate sequence $\{\xb_k\}$ generated by PGA-$\mathcal{B}$.
\begin{coro}\label{coro:rate of BFBS}
	Under Assumptions \ref{assum: compo} and \ref{assum: lipschitz}, the sequence $\{\xb_k\}$ generated by PGA-$\mathcal{B}$ for solving problem \textbf{(P2)} satisfies:
	\begin{flalign}
	&F(\xb_{k+1}) - F(\xb_{k}) \le - D_{\ell}(\xb_k, \xb_{k+1}),\label{eq: func_dec}\\
	& D_{\ell}(\xb^*, \xb_{k+1}) \le D_{\ell}(\xb^*, \xb_{k}), \quad \forall \xb^* \in \mathcal{X}^*, \label{eq: dist_dec}\\
	&F(\xb_{k}) - F(\ub) \le \frac{D_{\ell}(\ub,\xb_{0})}{k}, \quad\forall \ub \in C,\label{eq: global}
	\end{flalign}
\end{coro}

The result in \eqref{eq: func_dec} implies that the sequence of function value is non-increasing, and hence PGA-$\mathcal{B}$ is a descent method. Also, \eqref{eq: dist_dec} shows that the Bregman distance between $\xb_k$ and the optimal solution point $\xb^* \in \mathcal{X}^*$ is non-increasing. Moreover, \eqref{eq: global} with $\ub = \xb^* \in \mathcal{X}^*$ implies that the function value sequence $\{F(\xb_{k})\}$ converges to optimum at a rate $\mathcal{O}(1/k)$.

Similar results to those in Corollary \ref{coro:rate of BFBS} are established in \cite{bolte_BPGA}, but they are in terms of the Bregman distance $D_h$ (not $D_\ell$). Moreover, their analysis crucially depends on a symmetry coefficient $\alpha:= \inf_{\xb\ne \yb}\{\frac{D_h(\xb, \yb)}{D_h(\yb,\xb)}\} \in [0,1]$, which is avoided in our result through the unified point of view. Thus, our unification of PGA-$\mathcal{B}$ as PPA-$\mathcal{B}$ provides much simplicity of the analysis and avoids introducing the symmetry coefficient $\alpha$. Moreover, our global estimate in \eqref{eq: global} is tighter than the result in \cite[Theorem 1, (iv)]{bolte_BPGA}, since for all $\xb \in C$ and $\yb \in \mathrm{int} C$
\begin{align*}
D_{\ell}(\xb, \yb)\le \frac{1}{\gamma} D_{h}(\xb, \yb) \le \frac{2}{(1+\alpha)\gamma} D_h(\xb, \yb), \quad\forall \alpha \in [0,1].
\end{align*}

To further ensure the convergence of $\{\xb_k\}$ generated by PGA-$\mathcal{B}$ to a minimizer $\xb^*\in \mathcal{X}^*$, the following additional conditions on the Bregman distance $D_\ell$ are needed, and they are in parallel to the conditions introduced in \cite[Def 2.1, (iii)-(v)]{GPPA-Chen-Marc} to analyze the convergence of the sequence that generated by PPA-$\mathcal{B}$.
\begin{coro}
	Under Assumptions \ref{assum: compo} and \ref{assum: lipschitz}, the sequence $\{\xb_k\}$ generated by PGA-$\mathcal{B}$ for solving problem \textbf{(P2)} converges to some $\xb^* \in \mathcal{X}^*$ if the Bregman distance $D_{\ell}$ satisfies
	\begin{enumerate}[topsep=0pt,noitemsep]
		\item For every $\xb \in C$ and every $\alpha \in \RR$, the level set $\{\yb \in \mathrm{int}C ~|~ D_\ell(\xb, \yb) \le \alpha\}$ is bounded;
		\item If $\{\xb_k\} \in \mathrm{int}C$ and $\xb_k \to \xb^* \in C$, then $D_\ell(\xb^*, \xb_k) \to 0$;
		\item If $\{\xb_k\} \in \mathrm{int}C$ and $\xb^*\in C$ is such that $D_\ell(\xb^*, \xb_k) \to 0$, then $\xb_k \to \xb^*$.
	\end{enumerate}
\end{coro}
\begin{proof}
	The proof follows the argument in \cite[Theorem 3.4]{GPPA-Chen-Marc}.
\end{proof}

Next, we illustrate the PGA-$\mathcal{B}$ for an example, in which the objective function does not have a norm structure and satisfies Assumption \ref{assum: lipschitz}.
\begin{example}[Sparse Poisson Linear Inverse Problem]\label{ex: poisson}
	In many research areas, e.g., astronomy, electronic microscopy, one needs to solve inverse problems where the observations (e.g., photons, electrons) can be described by a Poisson process of the measurements of the underlying signal \cite{Bertero_2009}. Specifically, consider the following observation model:
	\begin{align}
	y_i = \mathrm{Poisson}(\inner{\ab_i}{\xb_0}),\quad i=1,\ldots,m, \label{eq: poisson}
	\end{align}
	where $\xb_0 \in \RR_+^n$ is the underlying signal to be recovered, $\{\ab_i\}_{i=1}^m$ are measurement vectors and $\{y_i\}_{i=1}^m$ correspond to observations of linear measurements under the Poisson noise.  
\end{example}
Denote $\yb := [y_1; \ldots; y_m]$, $A:=[\ab_1^\top; \ldots; \ab_m^\top]$ and assume that $\inner{\ab_i}{\xb_0}>0$ for all $i$, the goal is to recover $\xb_0$ given $\yb$ and $A$. Note that the corresponding negative Poisson log-likelihood function is given by
\begin{align*}
f(\xb) := D_p(\yb, A\xb), ~\text{where}~ p(\xb) = \sum_{j=1}^{n} x_j\log x_j - x_j.
\end{align*} 
Assume further that $\xb_{0}$ is sparse, we then add the regularizer $g(\xb) = \mu \|\xb\|_1, \mu>0$ to promote sparsity on the solution, and the optimization problem becomes
\begin{flalign}\label{eq:problem_poisson}
\textbf{(Q)}\quad \min_{\xb \in \RR_+^n} ~\{D_p(\yb, A\xb)+\mu \|\xb\|_1\}.
\end{flalign}
It has been shown in \cite[Lemma 7]{bolte_BPGA} that function $\frac{1}{\gamma} h - f$ is convex with $h(\xb) = -\sum_{i=1}^{n} \log x_i$ and $\gamma \le 1/\sum_{i=1}^{m} y_i$. Thus, PGA-$\mathcal{B}$ can be applied and the global estimate in \eqref{eq: global} holds. 
Moreover, with the specific choice of $h$, each iteration of PGA-$\mathcal{B}$ solves $m$ one-dimensional optimization problems, e.g., for all $j=1,\ldots,m$
\begin{align*}
&(\xb_{k+1})_j = \underset{x \in \RR_{++}}{\argmin} \left\{\mu x + \nabla_j f(\xb_k) x + \frac{1}{\gamma} \left(\frac{x}{(\xb_k)_j} - \log\frac{x}{(\xb_k)_j}\right)  \right\}.
\end{align*}
These one-dimensional subproblems have a closed-form solution characterized by: for all $j=1,\ldots,m$
\begin{align}
&(\xb_{k+1})_j = (\xb_{k})_j \left[1+\gamma (\xb_{k})_j \left(\mu+ \sum_{i=1}^{m} \left(A_{ij} - \frac{y_iA_{ij}}{\inner{\ab_i}{\xb_k}}\right)\right) \right]^{-1}. \label{eq: update_poisson}
\end{align}

\section{PGA-$\mathcal{B}$ with Line Search}\label{sec: line search}
The analysis in previous section requires to choose stepsize $\gamma < \bar{\gamma}$ in Assumption \ref{assum: lipschitz}, and $\bar{\gamma}$ is usually unknown a priori in practical applications. In particular,
it corresponds to the global Lipschitz parameter of $\nabla f$ when $h(\cdot) = \frac{1}{2}\|\cdot\|_2^2$. Next, we propose an adaptive version of PGA-$\mathcal{B}$ that searches for a proper stepsize in each step via the backtracking line search method. The algorithm is referred to as PGA-$\mathcal{B}$ with backtracking line search, and the details are summarized in Algorithm \ref{algo: pga-line}.
\begin{algorithm}
	\caption{PGA-$\mathcal{B}$ with backtracking line search} \label{algo: pga-line}
	Initialize $\gamma_0 > 0$ ,$\ell_0 = \frac{1}{\gamma_0}h - f$, $0<\beta<1$\;
	\For{$k=0, 1,\cdots$}
	{$\xb_{k+1} = \argmin_{\xb\in C} \{F(\xb) + D_{\ell_k}(\xb,\xb_k) \} $\;
		\While{$D_{\ell_k}(\xb_{k+1},\xb_k)<0$}
		{ Set $\gamma_k \leftarrow \beta \gamma_k$\;
			Repeat the $k$-th iteration with $\gamma_k$\;	
		}
	}
\end{algorithm}

We note that the above $D_{\ell_k}$ may not be a proper Bregman distance since $\ell_k$ may not be globally convex when $\gamma_k$ is large.
Intuitively, at each iteration we search for a small enough $\gamma_k$ that guarantees the non-negativity of the Bregman distance $D_{\ell_k}$ between the successive iterates $\xb_{k+1}$ and $\xb_k$. Intuitively, this implies that $\ell_k$ behaves like a convex function between these two successive iterates. Note that in the special case where $h(\cdot) = \frac{1}{2}\|\cdot\|_2^2$, the line search criterion $D_{\ell_k}(\xb_{k+1},\xb_k)\ge 0$ reduces to that of the PGA-$\mathcal{E}$ with line search in \cite{FISTA}. Next we characterize the convergence rate of PGA-$\mathcal{B}$ with line search.
\begin{thm}
	Under Assumptions \ref{assum: compo} and \ref{assum: lipschitz}, the sequence $\{\xb_k\}$ generated by PGA-$\mathcal{B}$ with backtracking line search satisfies: for all $k$ and all $\xb^* \in \mathcal{X}^*$
	\begin{align}
	F(\xb_{k}) - F^* \le \frac{D_{h}(\xb^*, \xb_0)}{\beta \bar{\gamma} k}.
	\end{align}
\end{thm}
\begin{proof}
	Since the line search method reduces $\gamma_k$ by a factor of $\beta$ whenever the Bregman distance is negative, we must have
	\begin{align}
	\beta \bar{\gamma} \le \gamma_k \le \bar{\gamma},\quad \forall k = 1, 2, \ldots.
	\end{align}
	At the $k$-th iteration, from PGA-$\mathcal{B}$ we have that
	\begin{align}
	F(\xb_{k+1}) + D_{\ell_k}(\xb_{k+1}, \xb_k) \le F(\xb_k),
	\end{align}
	which, combines with the line search criterion $D_{\ell_k}(\xb_{k+1}, \xb_k)\ge 0$, guarantees that $F(\xb_{k+1}) \le F(\xb_k)$, \ie, the method is a descent algorithm. Now by the convexity of $f$ and $g$, for any $\xb^* \in \mathcal{X}^*$ we have
	\begin{align}
	F(\xb^*) &\ge f(\xb_k) + \inner{\xb^* - \xb_k}{\nabla f(\xb_k)} + g(\xb_{k+1}) + \inner{\xb^* - \xb_{k+1}}{\partial g(\xb_{k+1})} \nonumber\\
	& = F(\xb_{k+1}) + \inner{\xb^* - \xb_{k+1}}{\partial F(\xb_{k+1})}+ D_f(\xb^*, \xb_{k+1}) - D_f(\xb^*, \xb_k). \label{eq: pga-line}
	\end{align}
	On the other hand, the optimality condition of the $(k+1)$-th iteration of PGA-$\mathcal{B}$ implies that
	\begin{align*}
	\nabla \ell_k (\xb_k) - \nabla \ell_k (\xb_{k+1}) \in \partial F(\xb_{k+1}),
	\end{align*}
	which together with the property in \eqref{eq:three point} further implies that
	\begin{align*}
	\inner{\xb^* - \xb_{k+1}}{\partial F(\xb_{k+1})} &= D_{\ell_k}(\xb^*, \xb_{k+1}) - D_{\ell_k}(\xb^*, \xb_{k}) + D_{\ell_k}(\xb_{k+1}, \xb_{k}) \\
	&\ge D_{\ell_k}(\xb^*, \xb_{k+1}) - D_{\ell_k}(\xb^*, \xb_{k}).
	\end{align*}
	Substituting the above inequality into \eqref{eq: pga-line}, we obtain that
	\begin{align*}
	0&\ge F(\xb^*) - F(\xb_{k+1}) \\
	&\ge \frac{1}{\gamma_k} \left[D_{h}(\xb^*, \xb_{k+1}) - D_{h}(\xb^*, \xb_{k}) \right] \\
	&\ge \frac{1}{\beta\bar{\gamma}} \left[D_{h}(\xb^*, \xb_{k+1}) - D_{h}(\xb^*, \xb_{k}) \right],
	\end{align*}
	where the last inequality follows from negativity and the fact that $\beta \bar{\gamma} \le \gamma_k$. Telescoping the above inequality from $0$ to $k-1$ and applying the fact that $F(\xb_{k+1}) \le F(\xb_k)$, we obtain that
	\begin{align*}
	k F(\xb^*) - k F(\xb_{k}) &\ge k F(\xb^*) - \sum_{l=1}^{k} F(\xb_{l}) \\
	&\ge \frac{1}{\beta\bar{\gamma}} \sum_{l=0}^{k-1}  \left[D_{h}(\xb^*, \xb_{l+1}) - D_{h}(\xb^*, \xb_{l}) \right] \\
	&\ge  - \frac{1}{\beta\bar{\gamma}} D_{h}(\xb^*, \xb_{0}).
	\end{align*}
	The result follows by rearranging the inequality.
\end{proof}

\section{Numerical Experiments}\label{sec: exp}

We now compare the practical performance between PGA-$\mathcal{B}$ and PGA-$\mathcal{B}$ with backtracking line search via numerical experiments. We consider applying these two algorithms to solve the sparse Poisson linear inverse problem $\textbf{(Q)}$ in Example \ref{ex: poisson}. 

Specifically, we randomly generate an underlying true signal $\xb_{0} \in \RR^{1000}$ with unit norm via the uniform distribution over the interval $[0,1]$. The sparsity of the signal is set to be either 1\% or 10\%. We consider two types of linear measurements, i.e., random measurements and convolution measurements. For the random measurements, we randomly generate 200 measurement vectors $\ab_1, \ldots, \ab_{200} \in \RR^{1000}$ via uniform distribution over the interval $[0,1]$. The observations $y_1, \ldots, y_{200}$ are obtained through the Poisson model in \eqref{eq: poisson}. For the convolution measurements, we first randomly generate the convolution kernel $\ub \in \RR^{20}$ via uniform distribution over the interval $[0,1]$. Then, we perform the linear convolution measurements $\ub * \xb_{0}$ and the observations are further obtained through the Poisson model in \eqref{eq: poisson}.
\begin{figure*}[htbp]
	\centering
	\begin{minipage}[t]{0.49\textwidth}
		\includegraphics[width=2.5in]{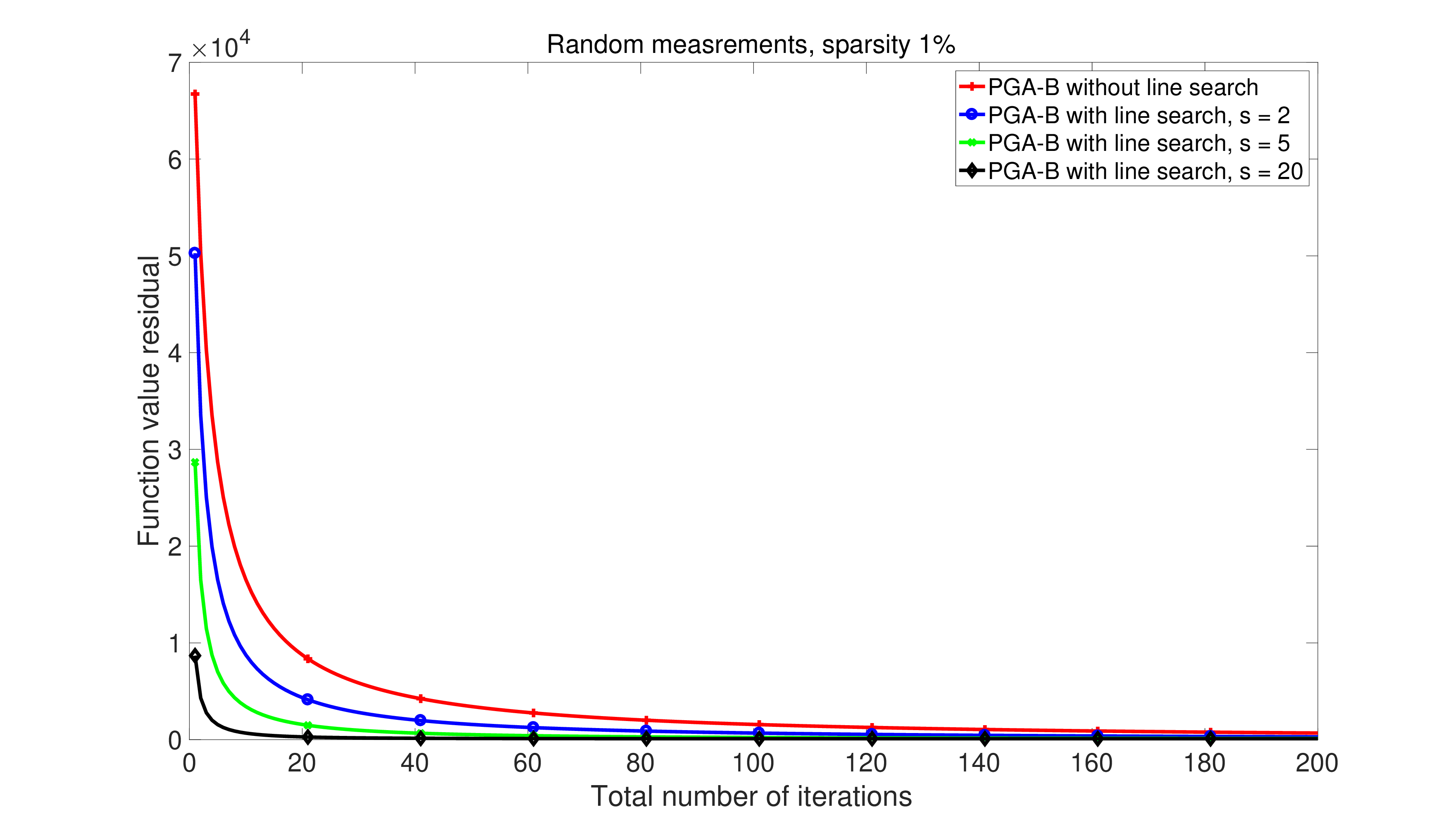}
	\end{minipage}
	\begin{minipage}[t]{0.49\textwidth}
		\includegraphics[width=2.5in]{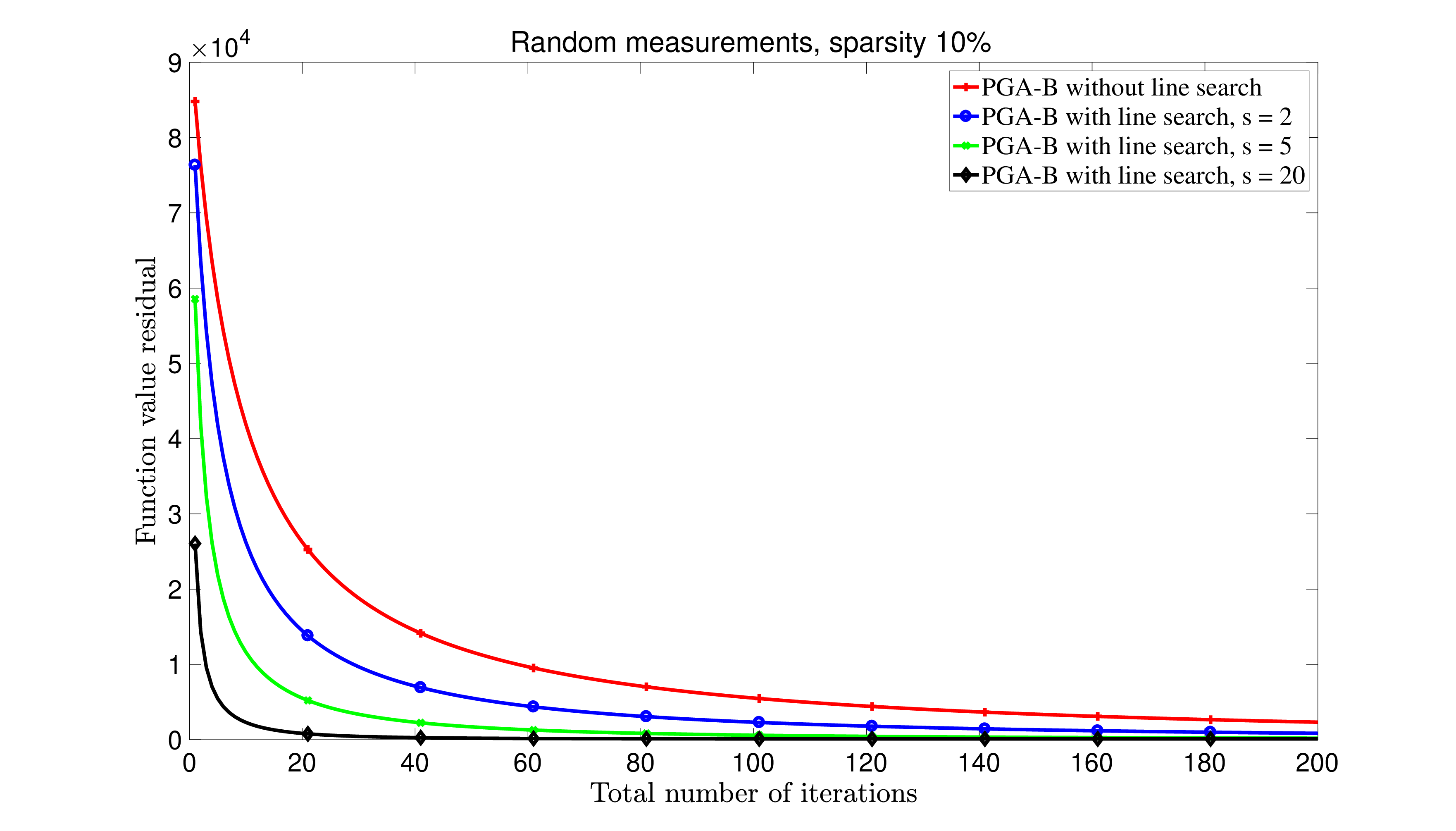}
	\end{minipage}
	\begin{minipage}[t]{0.49\textwidth}
		\includegraphics[width=2.5in]{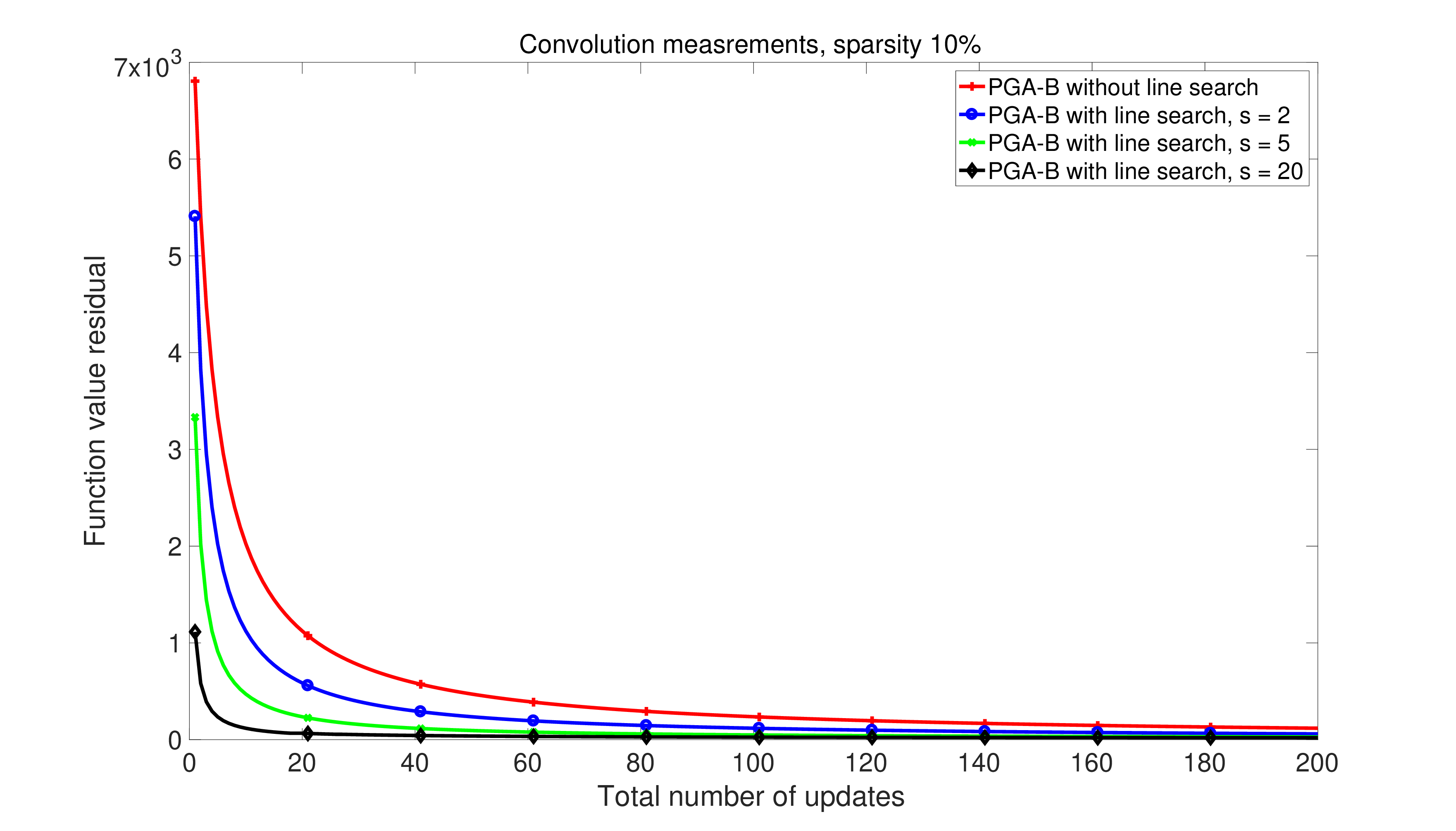}
	\end{minipage}
	\begin{minipage}[t]{0.49\textwidth}
		\includegraphics[width=2.5in]{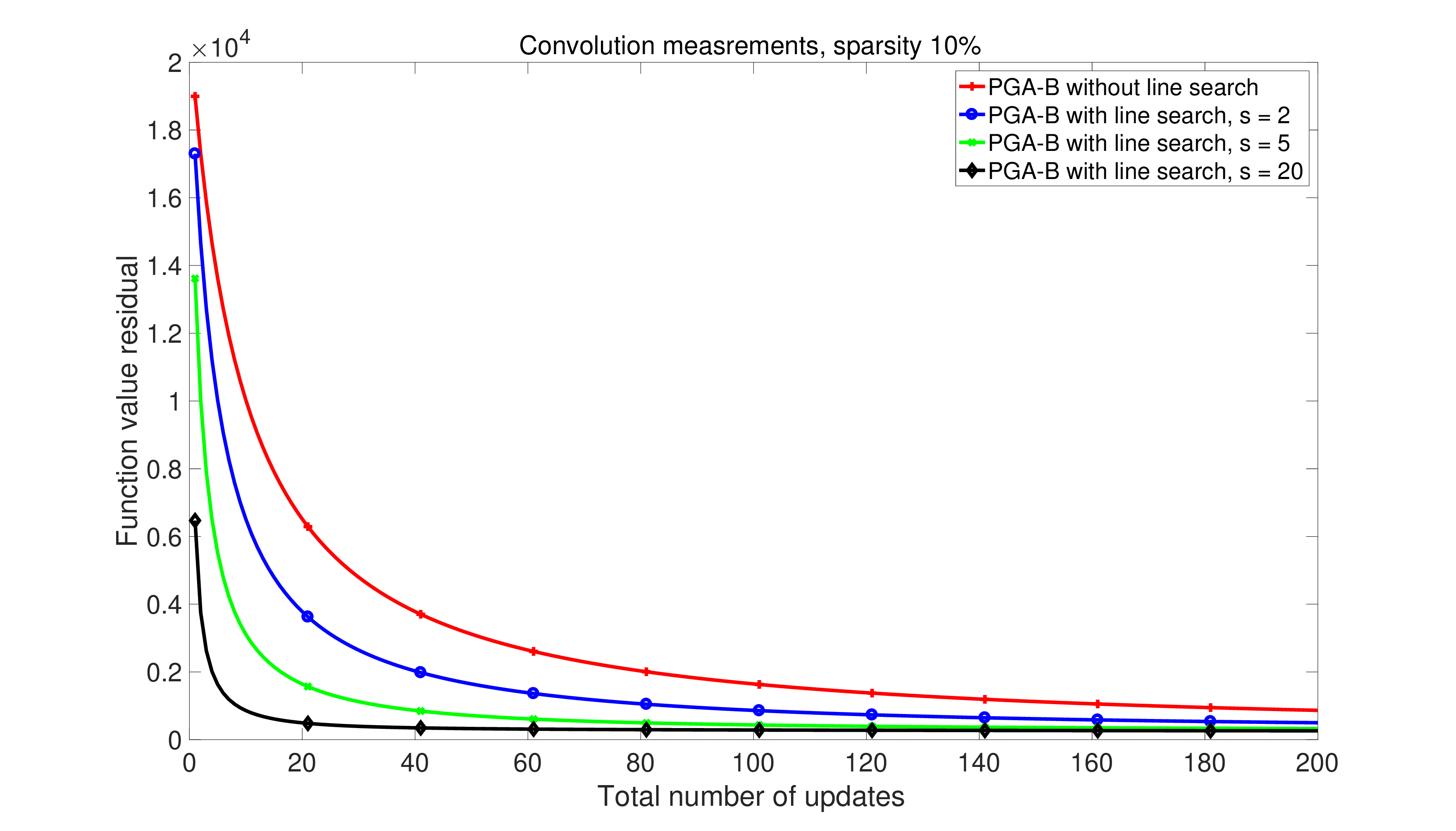}
	\end{minipage}
	\caption{Comparison of iteration complexity between PGA-$\mathcal{B}$ and PGA-$\mathcal{B}$ with line search under different settings.}\label{fig: 1}
\end{figure*}
We use the entropy function $h(\xb) = -\sum_{i=1}^{n} \log x_i$ for both PGA-$\mathcal{B}$ and PGA-$\mathcal{B}$ with backtracking line search, and the main update rules for both algorithms are given in \eqref{eq: update_poisson}. The stepsize for PGA-$\mathcal{B}$ is set to be $\gamma = 1/\sum_{i=1}^{m} y_i$. Moreover, for PGA-$\mathcal{B}$ with line search, we set the reducing factor for the stepsize to be $\beta = 0.8$ and consider different initializations of the stepsize: $\gamma_0 = s \gamma$ with $s = 2, 5, 20$. 

The residual of function value versus the number of  iterations for both algorithms are plotted in Figure \ref{fig: 1}, and we note that the repeated iterations of line search methods are also counted. It can be seen from the figure that 
PGA-$\mathcal{B}$ with backtracking line search converges much faster than vanilla PGA-$\mathcal{B}$ under different measurement settings and sparsity levels of the signal. This is because the line search method searches for a larger stepsize for the algorithm. Moreover, PGA-$\mathcal{B}$ benefits more from the line search method with a larger initialization stepsize $s\gamma$. This implies that the line search method quickly adapts the infeasible stepsize to a feasible one.

\section{Conclusion}\label{sec:conclude}
In this paper, we point out that PGA-$\mathcal{B}$ can be viewed as PPA-$\mathcal{B}$ with a special choice of Bregman distance. Consequently, the convergence analysis of PGA-$\mathcal{B}$ follows directly from that of PPA-$\mathcal{B}$. Moreover, this unified view point leads to a tighter convergence rate of the function value residual than existing results, and avoids involving the symmetry coefficient. Lastly, we provide a general line search variant of  PGA-$\mathcal{B}$ and characterize its convergence rate.

\bibliographystyle{apalike}      
\bibliography{./ref}   

\begin{thebibliography}{}

\bibitem[Auslender and Teboulle, 2006]{GPPA-Marc}
Auslender, A. and Teboulle, M. (2006).
\newblock Interior gradient and proximal methods for convex and conic
  optimization.
\newblock {\em SIAM Journal on Optimization}, 16(3):697--725.

\bibitem[Beck and Teboulle, 2003]{Tebolle_mirror}
Beck, A. and Teboulle, M. (2003).
\newblock Mirror descent and nonlinear projected subgradient methods for convex
  optimization.
\newblock {\em Operations Research Letters}, 31(3).

\bibitem[Beck and Teboulle, 2009]{FISTA}
Beck, A. and Teboulle, M. (2009).
\newblock A fast iterative shrinkage-thresholding algorithm for linear inverse
  problems.
\newblock {\em SIAM J. Img. Sci.}, 2(1):183--202.

\bibitem[Beck and Teboulle, 2012]{smooth-Marc}
Beck, A. and Teboulle, M. (2012).
\newblock Smoothing and first order methods: A unified framework.
\newblock {\em SIAM Journal on Optimization}, 22(2):557--580.

\bibitem[Bertero et~al., 2009]{Bertero_2009}
Bertero, M., Boccacci, P., DesiderA , G., and Vicidomini, G. (2009).
\newblock Image deblurring with poisson data: from cells to galaxies.
\newblock {\em Inverse Problems}, 25(12):123006.

\bibitem[Bolte et~al., 2016]{bolte_BPGA}
Bolte, J., Bauschke, H., and Teboulle, M. (2016).
\newblock A descent lemma beyond {L}ipschitz gradient continuity: first-order
  methods revisited and applications.
\newblock {\em Mathematics of Operations Research}.

\bibitem[Boyd et~al., 2011]{distributed-admm}
Boyd, S., Parikh, N., Chu, E., Peleato, B., and Eckstein, J. (2011).
\newblock Distributed optimization and statistical learning via the alternating
  direction method of multipliers.
\newblock {\em Found. Trends Mach. Learn.}, 3(1):1--122.

\bibitem[Bregman, 1967]{Bregman1967}
Bregman, L.~M. (1967).
\newblock The relaxation method of finding the common point of convex sets and
  its application to the solution of problems in convex programming.
\newblock {\em USSR Computational Mathematics and Mathematical Physics},
  7(3):200 -- 217.

\bibitem[Censor and Zenios, 1992]{Censor-GPPA}
Censor, Y. and Zenios, S.~A. (1992).
\newblock Proximal minimization algorithm with d-functions.
\newblock {\em J. Optim. Theory Appl.}, 73(3):451--464.

\bibitem[Chen and Teboulle, 1993]{GPPA-Chen-Marc}
Chen, G. and Teboulle, M. (1993).
\newblock Convergence analysis of a proximal-like minimization algorithm using
  {B}regman functions.
\newblock {\em SIAM Journal on Optimization}, 3(3):538--543.

\bibitem[De~Pierro and Iusem, 1986]{Pierro1986}
De~Pierro, A.~R. and Iusem, A.~N. (1986).
\newblock A relaxed version of bregman's method for convex programming.
\newblock {\em Journal of Optimization Theory and Applications},
  51(3):421--440.

\bibitem[Eckstein, 1993]{BPPA-Eckstein}
Eckstein, J. (1993).
\newblock Nonlinear proximal point algorithms using {B}regman functions, with
  applications to convex programming.
\newblock {\em Mathematics of Operations Research}, 18(1):202--226.

\bibitem[Eckstein and Bertsekas, 1992]{Eckstein-DR}
Eckstein, J. and Bertsekas, D.~P. (1992).
\newblock On the {D}ouglas-{R}achford splitting method and the proximal point
  algorithm for maximal monotone operators.
\newblock {\em Mathematical Programming}, 55:293--318.

\bibitem[Goldstein, 1964]{goldstein-PGA}
Goldstein, A.~A. (1964).
\newblock Convex programming in {H}ilbert space.
\newblock {\em Bulletin of the American Mathematical Society}, 70(5):709--710.

\bibitem[G\"{u}ler, 1991]{Guler-PPA}
G\"{u}ler, O. (1991).
\newblock On the convergence of the proximal point algorithm for convex
  minimization.
\newblock {\em SIAM J. Control Optim.}, 29(2):403--419.

\bibitem[Lions and Mercier, 1979]{splitting-algo}
Lions, P. and Mercier, B. (1979).
\newblock Splitting algorithms for the sum of two nonlinear operators.
\newblock {\em SIAM Journal on Numerical Analysis}, 16(6):964--979.

\bibitem[Martinet, 1970]{Martinet-PPA}
Martinet, B. (1970).
\newblock Br$\grave{\text{e}}$ve communication. r$\acute{\text{e}}$gularisation
  d'in$\acute{\text{e}}$quations variationnelles par approximations
  successives.
\newblock {\em ESAIM: Mathematical Modelling and Numerical Analysis -
  Mod$\acute{\text{e}}$lisation Math$\acute{\text{e}}$matique et Analyse
  Num$\acute{\text{e}}$rique}, 4(R3):154--158.

\bibitem[Micchelli et~al., 2011]{Micchelli-Shen-Xu:IP-11}
Micchelli, C.~A., Shen, L., and Xu, Y. (2011).
\newblock Proximity algorithms for image models: Denoising.
\newblock {\em Inverse Problems}, 27:045009(30pp).

\bibitem[Nesterov, 2005]{Nesterov-smooth}
Nesterov, Y. (2005).
\newblock Smooth minimization of non-smooth functions.
\newblock {\em Math. Program.}, 103(1):127--152.

\bibitem[Recht et~al., 2010]{low-rank}
Recht, B., Fazel, M., and Parrilo, P.~A. (2010).
\newblock Guaranteed minimum-rank solutions of linear matrix equations via
  nuclear norm minimization.
\newblock {\em SIAM Rev.}, 52(3):471--501.

\bibitem[Rockafellar, 1976]{RT-PPA}
Rockafellar, R. (1976).
\newblock Monotone operators and the proximal point algorithm.
\newblock {\em SIAM Journal on Control and Optimization}, 14(5):877--898.

\bibitem[Teboulle, 1992]{Entro-Marc}
Teboulle, M. (1992).
\newblock Entropic proximal mappings with applications to nonlinear
  programming.
\newblock {\em Mathematics of Operations Research}, 17(3):670--690.

\bibitem[Teboulle, 1997]{Marc-proximal-like}
Teboulle, M. (1997).
\newblock Convergence of proximal-like algorithms.
\newblock {\em SIAM Journal on Optimization}, 7(4):1069--1083.

\bibitem[Tseng, 2010]{Tseng-BPGA}
Tseng, P. (2010).
\newblock Approximation accuracy, gradient methods, and error bound for
  structured convex optimization.
\newblock {\em Mathematical Programming}, 125(2):263--295.

\bibitem[Zhou et~al., 2015]{Yi2015}
Zhou, Y., Liang, Y., and Shen, L. (2015).
\newblock A new perspective of proximal gradient algorithms.
\newblock {\em arXiv: 1503.05601}.

\end{thebibliography}

\end{document}